\documentclass{article}
\usepackage{amsmath,amssymb,amsthm,enumerate,float,graphicx,physics,tikz,tikz-cd,todonotes}

\usepackage{geometry}

\usepackage{hyperref}

\newtheorem{theorem}{Theorem}
\newtheorem{lemma}[theorem]{Lemma}
\newtheorem{corollary}[theorem]{Corollary}

\newtheorem{example}[theorem]{Example}

\theoremstyle{definition}
\newtheorem{definition}[theorem]{Definition}
\newtheorem{remark}[theorem]{Remark}

\newcommand{\PG}{\mathrm{PG}}

\newcommand{\FF}{\mathbb{F}}
\renewcommand{\a}{\mathbf{a}}
\renewcommand{\b}{\mathbf{b}}
\renewcommand{\c}{\mathbf{c}}
\newcommand{\e}{\mathbf{e}}
\newcommand{\f}{\mathbf{f}}
\newcommand{\p}{\mathbf{p}}
\newcommand{\x}{\mathbf{x}}

\renewcommand{\u}{\mathbf{u}}
\renewcommand{\v}{\mathbf{v}}

\newcommand{\0}{\mathbf{0}}

\usepackage{xcolor}

\title{Stabiliser codes over prime fields\\ as weighted graphs and as quantum sets of lines}
\author{Simeon Ball\thanks{Supported by the Spanish Ministry of Science, Innovation and Universities grant PID2023-147202NB-I00.}\qquad Gülizar Günay\thanks{This research was funded in whole or in part by the Austrian Science Fund (FWF) [10.55776/ESP4162324].}\qquad Robin Simoens\thanks{Supported by the Research Foundation Flanders (FWO) through the grant 11PG724N.}\qquad Ricard Vilar}
\date{}

\begin{document}

\maketitle

\begin{abstract}
\noindent
    We further develop the theory of quantum stabiliser codes and its links with graph theory and finite geometry.
    We prove that an \([\![n,k,d]\!]_p\) stabiliser code over a finite field with \(p\) elements, \(p\) prime, gives rise to a graph on \(n-k\) vertices together with two disjoint independent sets of size \(k\). Its edges are weighted with elements of a finite field of size \(p\). If \(k=0\), such codes are known as graph states.
    Second, we view stabiliser codes as so-called quantum sets of lines in projective space, allowing us to determined their minimum distance geometrically. We describe two geometrical operations to obtain new quantum sets of lines from given quantum sets of lines. As an application, we list the number of stabiliser codes up to equivalence, for given small parameters.
\end{abstract}

\section{Introduction}

Stabiliser codes were introduced by Calderbank, Rains, Shor and Sloane \cite{CRSS1998} and Gottesman \cite{gottesman96} independently. They are a special type of quantum error-correcting codes.
A \emph{stabiliser code} is a subspace \(Q\) of a Hilbert space \(\left(\mathbb{C}^q\right)^{\otimes n}\) that is the intersection of the eigenspaces of eigenvalue one of an abelian subgroup \(S\) of \(\mathcal{E}^{\otimes n}\), where \(\mathcal{E}\) is a nice error basis. 
These codes are widely used in quantum error-correction, partly because it is often easier to work with the subgroup \(S\), rather than with a basis for the subspace \(Q\). In this article, we restrict ourselves to the case where \(q\) is prime. In this case, the stabiliser formalism allows us to interpret the abelian subgroup \(S\) as a classical additive code that is contained in its symplectic dual, see \cite{AK2001,BCH2023,CRSS1998,KKKS2006}.

One-dimensional binary stabiliser codes are called graphical states or \emph{graph states} since they can be represented by a simple graph.
Graph states have been well studied, see for example \cite{AB2006,Elliot2008,HEB2004,KLS2025}. There have been many proposals to use graphs to represent stabiliser codes of general dimension, including \cite{HK2022,MB2022,KLS2025,YCO2007}.

In this article, we prove that an \([\![n,k,d]\!]_p\) stabiliser code over a prime field can be interpreted as a weighted graph on \(n+k\) vertices containing two independent sets of size \(k\) that induce a (non-weighted) matching. As far as we are aware, this equivalence has not been previously noted.

We go on to consider \emph{quantum sets of lines}, sets of lines in projective space, associated to such a stabiliser code, from which it is possible to give a geometrical description of the minimum distance $d$. These ideas are inspired by the work of Glynn, Gulliver, Maks and Gupta \cite{GGMG} and build on the generalisations of Ball, Centelles and Huber \cite{BCH2023}. See also the work of Bierbrauer, Faina, Giulietti, Marcugini and Pambianco
\cite{geometry}.

As an application, we present some tables on the number of stabiliser codes of dimension \(k=1\) and local dimension \(q=2\) (qubit) and \(q=3\) (qutrit), expanding on the known classifications for \(k=0\) and \(q=2\) by Danielsen and Parker \cite{classification} and \(k=0\) and \(q=3\) by Danielsen \cite{classification1}.

\section{Preliminaries}


\subsection{Stabiliser formalism}

Let \(\FF_p\) denote the finite field of order \(p\), where \(p\) is prime.
Consider the standard \emph{symplectic inner product} on the vector space \(\FF_p^{2n}\) given by
\[\langle\u,\v\rangle=\sum_{i=1}^nu_iv_{n+i}-v_iu_{n+i}.\]
Given a subspace \(U\subseteq\FF_p^{2n}\), let \(U^\perp=\{\v\in\FF_p^{2n}\colon\,\langle\u,\v\rangle=0\text{ for all }\u\in U\}\). A subspace \(U\) is called \emph{totally isotropic} if \(U\subseteq U^\perp\). The \emph{symplectic weight} of a vector \(\v\in\mathbb{F}_p^{2n}\) is defined as
\[\mathrm{swt}(\v)=\#\left\{i\in\{1,\dots,n\}\colon\,(v_i,v_{n+i})\neq(0,0)\right\}.\]

\begin{definition}
    An \emph{\([\![n,k]\!]_p\) stabiliser code} is a totally isotropic subspace \(S\) of \(\FF_p^{2n}\) of dimension \(n-k\). The \emph{minimum distance} of \(S\) is
    \[d(S)=\begin{cases}\min\left\{\mathrm{swt}(\v)\colon\,\v\in S^\perp\setminus S\right\}&\text{if } k>0,\\\min\left\{\mathrm{swt}(\v)\colon\,\v\in S\setminus\{\0\}\right\}&\text{if }k=0.\end{cases}\]
    If the minimum distance \(d\) is specified, we call \(S\) an \([\![n,k,d]\!]_p\) code.
\end{definition}

\begin{definition}
    A \emph{stabiliser matrix} of an \([\![n,k]\!]_p\) stabiliser code \(S\) is an \((n-k)\times2n\) matrix whose rows form a basis of \(S\).
    A \emph{parity check matrix}\footnote{We use the term ``parity check matrix'' in analogy with the term for classical codes, when considering the stabiliser matrix as a generator matrix. However, one could also regard the stabiliser matrix as the parity check matrix. After all, the minimum weight depends mostly on the minimum symplectic weight of \(H\).} of an \([\![n,k]\!]_p\) stabiliser code \(S\) is an \((n+k)\times2n\) matrix whose rows form a basis of \(S^\perp\).
\end{definition}

We denote a stabiliser matrix by \(G\) and a parity check matrix by \(H\).

\subsection{Equivalent stabiliser codes}

\begin{definition}\label{def:equivalent}
    Two stabiliser codes are \emph{(local Clifford) equivalent} if they can be obtained from one another by:
    \begin{enumerate}
        \item Permutations of the \(n\) pairs of positions.
        \item For each pair, a \emph{symplectic transformation} on that pair. That is, a change of basis of \(\FF_p^2\) that preserves the symplectic inner product. 
    \end{enumerate}
\end{definition}
The maps yielding equivalence are the elements of the \emph{local Clifford group} \(\mathrm{Sp}(2,p)^n\rtimes\mathrm{Sym}(n)\).
For example, for a given \(i\in\{1,\dots,n\}\), we can replace \((v_i,v_{n+i})\) by \((v_{n+i},-v_i)\) within every vector \(\v\).
Another possibility is to replace \((v_i,v_{n+i})\) by \((v_i+v_{n+i},v_{n+i})\) within every vector \(\v\). These two operations are in fact enough to generate the whole group of symplectic transformations, since
\(\mathrm{Sp}(2,p)=\mathrm{SL}(2,p)\) is generated as a multiplicative group by all products of the matrices \[\begin{pmatrix}0&-1\\1&0\end{pmatrix}\text{ and } \begin{pmatrix}1&0\\1&1\end{pmatrix}.\]

Equivalent codes have the same parameters.

\section{Stabiliser codes as graphs}

A particular type of one-dimensional (\(k=0\)) stabiliser codes are \emph{graph states}. Given a graph on \(n\) vertices, we can construct from it a stabiliser matrix that has the identity matrix on the left and the adjacency matrix \(A\) of the graph on the right:
\[G=\left(\begin{array}{c||c}I_n&A\end{array} \right).\]
Bahramgiri and Beigi \cite{BB2006} proved that every \([\![n,0]\!]_p\) stabiliser code, \(p\) prime, is equivalent to a graph state. This was proved already by Schlingemann \cite{graphstate} for the binary case.
We generalise this result to multidimensional (\(k\geq1\)) stabiliser codes.

\subsection{Multidimensional graph states}

\begin{theorem}\label{thm:ABC}
Consider an \([\![n,k,d]\!]_p\) stabiliser code and let \(r=n-k\). Up to equivalence, its stabiliser matrix has the form
\[
G=\left(\begin{array}{c|c||c|c} I_r & B & A+BC^T & C \end{array} \right)
\]
where \(A\) is a symmetric \(r\times r\) matrix with zeroes on the diagonal and  \(B\) and \(C\) are \(r\times k\) matrices.
\end{theorem}

\begin{proof}
After row reduction and permutations of pairs of columns, we may assume that the code has stabiliser matrix
\[
G=\left(\begin{array}{c|c||c|c} I_{r'}&B&D&C
\\ \hline O & O & E & O\\
\end{array} \right),
\]
for some \(r'\leq r\). Here, the part below \(C\) may be assumed zero, since otherwise we can swap the columns containing \(B\) and \(C\) and increase the rank \(r'\) of the left part.
The rowspace of \(G\) is symplectic self-orthogonal if and only if
\(G\Lambda G^T=O\) where \(\Lambda=\left(\begin{array}{cc}O&I_n\\-I_n&O\end{array}\right)\). Thus,
\[O=G\Lambda G^T=\left(\begin{array}{c|c||c|c} I_{r'}&B&D&C
\\ \hline O & O & E & O\\
\end{array} \right)\left(\begin{array}{c|c} D^T&E^T\\
\hline
C^T&O\\ \hline\hline
-I_{r'}&O\\\hline
-B^T&O
\end{array} \right)=\left(\begin{array}{c|c}
D^T+BC^T-D-CB^T&E^T\\
\hline
-E&O\end{array}\right).
\]
We deduce that \(E=O\), and hence \(r'=r\), and that \(A=A^T\), where \(A=D-BC^T\). Finally, we can replace any non-zero elements on the diagonal of the matrix \(A\) by zero, by adding the \(i\)th column to the \((n+i)\)th column if needed (preserving the symplectic inner product of rows).
\end{proof}

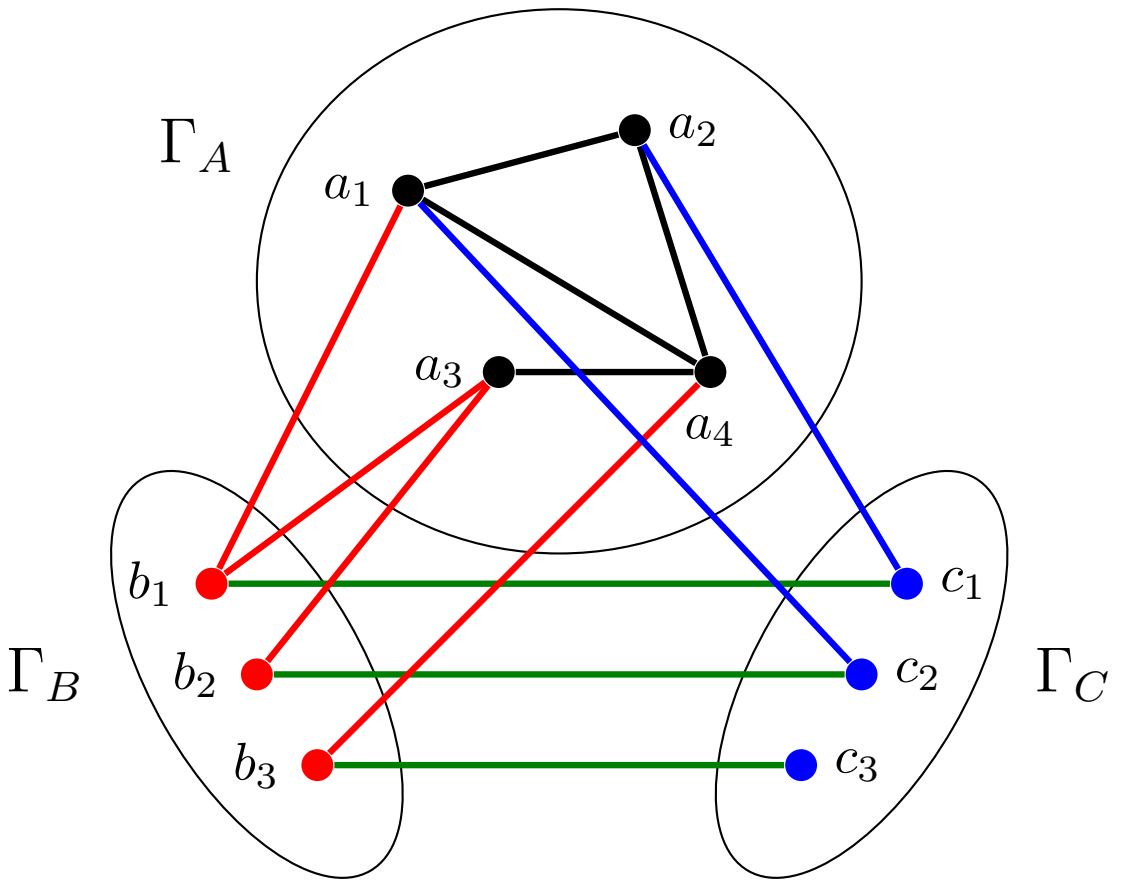
\begin{figure}[H]
    \centering
    \newcommand{\radius}{1}
\begin{tikzpicture}
    \draw[rotate around={60:(2,0)}] (2,0) ellipse (1.5 and .7);
    \draw[rotate around={-60:(-2,0)}] (-2,0) ellipse (1.5 and .7);
    \draw (0,2.6) ellipse (2 and 1.8);
    \draw (3.4,0) node {\large\(\Gamma_C\)};
    \draw (-3.4,0) node {\large\(\Gamma_B\)};
    \draw (-2.4,3.5) node {\large\(\Gamma_A\)};
	\path[every node/.append style={circle, fill=black, minimum size=6pt, label distance=2pt, inner sep=0pt}]
    (.5,3.6) node[label={0:\(a_2\)}] (a4) {}
    (-1,3.2) node[label={180:\(a_1\)}] (a2) {}
    (1,2) node[label={-90:\(a_4\)}] (a3) {}
    (-.4,2) node[label={180:\(a_3\)}] (a1) {}
    (-2.3,.6) node[red,label={180:\(b_1\)}] (b1) {}
    (-2,0) node[red,label={180:\(b_2\)}] (b2) {}
    (-1.6,-.6) node[red,label={180:\(b_3\)}] (b3) {}
    (2.3,.6) node[blue,label={0:\(c_1\)}] (c1) {}
    (2,0) node[blue,label={0:\(c_2\)}] (c2) {}
    (1.6,-.6) node[blue,label={0:\(c_3\)}] (c3) {};
    \draw[line width=1.2] (a2) edge (a4) edge (a3) edge[red] (b1)
    (a3) edge (a1) edge (a4)
    (a1) edge[red] (b1)
    (c1) edge[blue] (a4) edge[green!50!black] (b1)
    (b2) edge[green!50!black] (c2) edge[red] (a1)
    (b3) edge[green!50!black] (c3) edge[red] (a3)
    (c2) edge[blue] (a2);
\end{tikzpicture}
    \caption{A multidimensional graph state is represented by a graph on \(n+k\) vertices containing two independent sets of size \(k\) with a matching in between. In this figure, \(n=7\), \(k=3\) and \(p=2\) (the graph is not weighted).}
    \label{fig:graphstate}
\end{figure}

More explicitly, we can always put the stabiliser matrix in the form

\[\left(\begin{array}{@{}ccc|ccc||ccc|ccc@{}}
1&\cdots&0&b_{11}&\cdots&b_{1k}&a_{11}+\sum_{i=1}^kb_{1i}c_{1i}&\cdots&a_{1r}+\sum_{i=1}^kb_{1i}c_{ri}&c_{11}&\cdots&c_{1k}\\
\vdots&\ddots&\vdots&\vdots&&\vdots&\vdots&\ddots&\vdots&\vdots&&\vdots\\
0&\cdots&1&b_{r1}&\cdots&b_{rk}&a_{r1}+\sum_{i=1}^kb_{ri}c_{1i}&\cdots&a_{rr}+\sum_{i=1}^kb_{ri}c_{ri}&c_{r1}&\cdots&c_{rk}
\end{array}\right)\]
where \(a_{ii}=0\) for all \(i\in\{1,\dots,r\}\).
The matrix \(A\) gives rise to an (undirected, loopless) weighted graph \(\Gamma_A\) on \(n-k\) vertices with weights in \(\mathbb{F}_p\). The matrices \(B\) and \(C\) can be seen as two empty graphs \(\Gamma_B\) and \(\Gamma_C\) on \(k\) vertices each, with a matching between them (this matching has no weights) and weighted edges between \(\Gamma_A\) and \(\Gamma_B\) (resp.\ \(\Gamma_C\)), indicating the entries of \(B\) (resp.\ \(C\)). Let \(\Gamma\) be the union (not disjoint) of \(\Gamma_A\), \(\Gamma_B\) and \(\Gamma_C\). An example is given in Figure~\ref{fig:graphstate}.
Label its vertices as follows:
\begin{align*}
    V_A&=\{a_1,\dots,a_r\},\\
    V_B&=\{b_1,\dots,b_k\},\\
    V_C&=\{c_1,\dots,c_k\}.
\end{align*}
Beware of the difference in notation with weights of the edges between the vertices. For example, \(a_{12}\) is the weight of the edge between the vertices \(a_1\) and \(a_2\).

\begin{example}
Consider the \([\![5,1,3]\!]_2\) stabiliser code generated by
\[
G=\left(\begin{array}{cccc|c||cccc|c} 
1&0&0&1&1&1&0&1&0&1\\ 0&1&0&0&1&1&0&0&1&1\\ 0&0&1&1&0&1&0&0&1&1\\ 0&0&0&0&0&0&1&1&1&1
 \end{array} \right).
\]
Swapping column 4 and column 9 and using row operations, yields an equivalent code with generator matrix
\[
\left(\begin{array}{cccc|c||cccc|c} 
1&0&0&0&1&1&0&1&1&1 \\ 0&1&0&0&1&1&1&1&0&0 \\ 0&0&1&0&0&1&1&1&1&0 \\ 0&0&0&1&0&0&1&1&0&1
 \end{array} \right)
\]
which, by adding columns 2-3 to columns 7-8, is equivalent to
\[
\left(\begin{array}{cccc|c||cccc|c} 
1&0&0&0&1&1&0&1&1&1 \\ 0&1&0&0&1&1&0&1&0&0 \\ 0&0&1&0&0&1&1&0&1&0 \\ 0&0&0&1&0&0&1&1&0&1
 \end{array} \right).
\]
We deduce that
\[
A=\begin{pmatrix}
0&0&1&0\\ 
0&0&1&1\\
1&1&0&1\\
0&1&1&0
\end{pmatrix},
\,B=\begin{pmatrix}1\\1\\0\\0\\\end{pmatrix}\text{ and }C=\begin{pmatrix}1\\0\\0\\1\\\end{pmatrix}.
\]
The corresponding graph is given in Figure~\ref{graph1}.
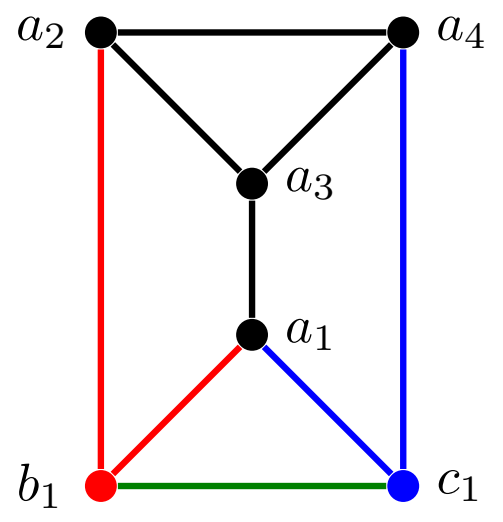
\begin{figure}[H]
\begin{center}
\newcommand{\radius}{1}
\begin{tikzpicture}
	\path[every node/.append style={circle, fill=black, minimum size=6pt, label distance=2pt, inner sep=0pt}]
    (2*\radius,3*\radius) node[label={0:\(a_4\)}] (a4) {}
    (0,3*\radius) node[label={180:\(a_2\)}] (a2) {}
    (\radius,2*\radius) node[label={0:\(a_3\)}] (a3) {}
    (\radius,\radius) node[label={0:\(a_1\)}] (a1) {}
    (0,0) node[red,label={180:\(b_1\)}] (b1) {}
    (2*\radius,0) node[blue,label={0:\(c_1\)}] (c1) {};
    \draw[line width=1.2] (a2) edge (a4) edge (a3) edge[red] (b1)
    (a3) edge (a1) edge (a4)
    (a1) edge[red] (b1) edge[blue] (c1)
    (c1) edge[blue] (a4) edge[green!50!black] (b1);
\end{tikzpicture}
\caption{A graph describing a \([\![5,1,3]\!]_2\) stabiliser code.}
\label{graph1}
\end{center}
\end{figure}
\end{example}

Vice versa, any graph on \(n+k\) vertices and two disjoint copies of the empty graph on $k$ vertices, with a matching between them and with edges joined to the subset of $n-k$ vertices, gives rise to an \([\![n,k]\!]\) stabiliser code. We call such codes \emph{multidimensional graph states}. Theorem~\ref{thm:ABC} says that every stabiliser code is equivalent to a multidimensional graph state.

\subsection{Parity check matrix}

Given the standard form
\[G=\left(\begin{array}{c|c||c|c}I_{n-k}&B&A+BC^T&C\end{array}\right)\]
from Theorem~\ref{thm:ABC}, we can determine a parity check matrix as follows.

\begin{lemma}\label{lem:pcm1}
    An \([\![n,k,d]\!]\) code with stabiliser matrix
    \[G=\left(\begin{array}{c|c||c|c}I_{n-k}&B&D&C\end{array}\right)\]
    has parity check matrix
    \[H=\left(\begin{array}{c|c||c|c}
    I_{n-k}&O&D^T&O\\
    \hline
    O & I_k & C^T & O \\
    \hline
    O & O & -B^T & I_k
    \end{array}\right).\]
\end{lemma}
\begin{proof}
    \(H\) has full rank and \(n+k\) rows.
    We check that its rows are symplectic orthogonal with those of \(G\). We have
    \[
    G\Lambda H^T=\left(\begin{array}{c|c||c|c}I_{n-k}&B&D&C\end{array}\right)\left( \begin{array}{c|c|c}D&C&-B\\ \hline O&O&I_k\\ \hline\hline-I_{n-k}&O&O\\ \hline O&-I_k&O \end{array} \right)=\left(\begin{array}{c|c|c}
    D-D& C-C & -B+B
    \end{array}\right)=O.\]
\end{proof}

A parity check matrix is not unique. The following construction has the advantage that the parity check matrix contains \(G\).

\begin{lemma}\label{lem:pcm2}
    An \([\![n,k,d]\!]\) code with stabiliser matrix
    \[G=\left(\begin{array}{c|c||c|c}I_{n-k}&B&D&C\end{array}\right)\]
    has parity check matrix
    \[H=\left(\begin{array}{c|c||c|c}
    I_{n-k}&B&D&C\\
    \hline
    O & I_k & C^T & O \\
    \hline
    O & O & -B^T & I_k
    \end{array}\right).\]
\end{lemma}
\begin{proof}
    \(H\) has full rank and \(n+k\) rows. The first \(n-k\) rows form the matrix \(G\) itself, which is symplectic self-orthogonal. The remaining \(2k\) rows are the same last \(2k\) rows as in Lemma~\ref{lem:pcm1}.
\end{proof}

\subsection{The Calderbank-Shor-Steane construction}

The CSS construction of stabiliser codes is a construction of stabiliser codes introduced by Calderbank and Shor \cite{CS1996} and Steane \cite{Steane1996} using two classical linear codes where the dual of one is contained in the other. The stabiliser matrix \(G\) for this construction has the form
\[G=\left(\begin{array}{c||c}G_1&O\\\hline O&G_2 \end{array} \right)\]
where \(G_1\) and \(G_2\) have full rank and \(G_1G_2^T=0\).

We proved in Theorem~\ref{thm:ABC} that an \([\![n,k]\!]\) stabiliser code is equivalent to a graph $\Gamma$ on $n+k$ vertices which has two copies of the empty graph on $k$ vertices and an induced subgraph on $n-k$ vertices. In the following theorem, we prove that the graph $\Gamma$ for the CSS construction is bipartite.

\begin{theorem} \label{ABCthmCSS}
A stabiliser matrix is obtained from a CSS construction if and only if, up to equivalence, it is of the form
\[
G=\left(\begin{array}{c|c||c|c} I_{n-k} & B & A+BC^T & C \end{array} \right)
\]
where
\[A=
\left(\begin{array}{c|c}O& A'\,^T \\ \hline A' & O\end{array} \right),
\,B=\left(\begin{array}{c}B'\\\hline O\end{array}\right)
\text{ and }
C=\left(\begin{array}{c}O\\\hline C'\end{array}\right).
\]
\end{theorem}
\begin{proof}
Let \(G=\left(\begin{array}{c||c}G_1&O\\\hline O&G_2 \end{array} \right)\) be a stabiliser matrix constructed using the CSS construction. We can use row operations and a re-ordering of the columns to put it in the form
\[G\sim\left(\begin{array}{c|c|c||c|c|c}I_{r_1}&G_{12}&G_{13}&O&O&O\\\hline O&O&O&G_{21}&I_{r_2}&G_{23}\\\hline O&O&O&E&O&O\end{array}\right)\]
where \(r_1\) is the rank of \(G_1\) and \(r_2+\mathrm{rank}(E)\) is the rank of \(G_2\).
Since the symplectic inner product of the $i$-th row ($i \in \{1,\dots,r_1\}$) and the $j$-th row (\(j \in \{r_1+r_2+1,\dots,n-k\}\)) is zero, we deduce that \(E=O\) and hence \(r_1+r_2=n-k\).
Swapping (and adding a minus sign on one side) the $i$th and \((n+i)\)th column for all \(i\in \{r_1+1,\dots,r_1+r_2\}\), we obtain
\[
G\sim\left(\begin{array}{c|c|c||c|c|c}I_{r_1}&O&G_{13}&O&-G_{12}&O\\\hline O&I_{r_2}&O&G_{21}&O&G_{23}\end{array} \right).
\]
Comparing this matrix with the matrix from Theorem~\ref{thm:ABC}, we get
\[A=\left(\begin{array}{c|c}O&-G_{12}-G_{13}G_{23}^T\\\hline G_{21}&O\end{array}\right),
\,B=\left(\begin{array}{c}G_{13}\\\hline O\end{array}\right)
\text{ and }
C=\left(\begin{array}{c}O\\\hline G_{23}\end{array}\right),
\]
and \(A\) is symmetric. By doing the inverse operations, a stabiliser matrix is of the given form can be obtained from the CSS construction with \(G_1=\left(\begin{array}{c|c|c}I_{r_1}&-A'\,^T-B'C'\,^T&B'\end{array}\right)\) and \(G_2=\left(\begin{array}{c|c|c}A'&I_{r_2}&C'\end{array}\right)\).
\end{proof}

As in the discussion after the proof of Theorem~\ref{thm:ABC}, we can interpret $A$ as describing an induced subgraph \(\Gamma_A\) of a graph $\Gamma$ on $n+k$ vertices. For the CSS construction, Theorem~\ref{ABCthmCSS} implies that the graph $\Gamma_A$ is bipartite with bipartitions of size $r_1$ and $r_2$. The graph $\Gamma$ is bipartite with bipartitions of size $r_1+k$ and $r_2+k$, where each bipartition contains an independent set on $k$ vertices.

\begin{example}
The Steane code is a $[[7,1,3]]_2$ stabiliser code, obtained via a CSS construction where
\[G_1=G_2=\begin{pmatrix}
    1&0&0&1&0&1&1\\0&1&0&1&1&1&0\\0&0&1&0&1&1&1
\end{pmatrix}\]
is the parity check matrix of the binary \([7,4,3]\) Hamming code. Its stabiliser matrix is
\[G=\left(\begin{array}{ccccccc||ccccccc} 
1 & 0 & 0 &1&0& 1 & 1 & 0 & 0 & 0 & 0 & 0 & 0 & 0 \\ 
0 & 1 & 0 & 1 & 1 & 1 & 0 & 0 & 0 & 0 & 0 & 0 & 0 & 0   \\ 
0 & 0 & 1 &0 & 1 &1 & 1 & 0 & 0 & 0 & 0 & 0 & 0 & 0\\
\hline
0 & 0 & 0 & 0 & 0 & 0 & 0 & 1 & 0 & 0 & 1 & 0 & 1 & 1  \\ 
0 & 0 & 0 & 0 & 0 & 0 & 0  & 0 & 1 & 0 & 1 & 1 & 1 & 0  \\ 
0 & 0 & 0 & 0 & 0 & 0 & 0  & 0 & 0 & 1 & 0 & 1 & 1 & 1 \\ 
 \end{array} \right).\]
By elementary row operations on the last three rows, we get
\[G\sim\left(\begin{array}{ccccccc||ccccccc} 
1 & 0 & 0 &1&0& 1 & 1 & 0 & 0 & 0 & 0 & 0 & 0 & 0 \\ 
0 & 1 & 0 & 1 & 1 & 1 & 0 & 0 & 0 & 0 & 0 & 0 & 0 & 0   \\ 
0 & 0 & 1 &0 & 1 &1 & 1 & 0 & 0 & 0 & 0 & 0 & 0 & 0\\
\hline
0 & 0 & 0 & 0 & 0 & 0 & 0 & 0 & 1 & 1 & 1 & 0 & 0 & 1\\ 
0 & 0 & 0 & 0 & 0 & 0 & 0 & 1 & 1 & 0 & 0 & 1 & 0 & 1\\ 
0 & 0 & 0 & 0 & 0 & 0 & 0 & 1 & 1 & 1 & 0 & 0 & 1 & 0\\ 
 \end{array} \right).
\]
Switching columns \(4\), \(5\) and \(6\) with columns \(11\), \(12\) and \(13\) respectively, we get
\[G\sim\left(\begin{array}{ccccccc||ccccccc} 
1 & 0 & 0 & 0 & 0 & 0 & 1 & 0 & 0 & 0 & 1 & 0 & 1 & 0\\ 
0 & 1 & 0 & 0 & 0 & 0 & 0 & 0 & 0 & 0 & 1 & 1 & 1 & 0\\ 
0 & 0 & 1 & 0 & 0 & 0 & 1 & 0 & 0 & 0 & 0 & 1 & 1 & 0\\ 
0 & 0 & 0 & 1 & 0 & 0 & 0 & 0 & 1 & 1 & 0 & 0 & 0 & 1\\ 
0 & 0 & 0 & 0 & 1 & 0 & 0 & 1 & 1 & 0 & 0 & 0 & 0 & 1\\ 
0 & 0 & 0 & 0 & 0 & 1 & 0 & 1 & 1 & 1 & 0 & 0 & 0 & 0 \\ 
\end{array} \right).
\]
Using the notation of Theorem~\ref{thm:ABC}, we have
\[A=\begin{pmatrix}
 0 & 0 & 0 & 0 & 1 & 1\\ 
 0 & 0 & 0 & 1 & 1 & 1\\ 
 0 & 0 & 0 & 1 & 0 & 1\\ 
 0 & 1 & 1 & 0 & 0 & 0\\ 
 1 & 1 & 0 & 0 & 0 & 0\\ 
 1 & 1 & 1 & 0 & 0 & 0\\ 
\end{pmatrix},
\,B=\begin{pmatrix}1\\0\\1\\0\\0\\0\\\end{pmatrix}\text{ and }C=\begin{pmatrix}0\\0\\0\\1\\1\\0\\\end{pmatrix}.
\]
Therefore, the Steane code is given by the bipartite graph in Figure~\ref{graphsteane}.
\begin{figure}[H]
\begin{center}
\newcommand{\radius}{1}
\begin{tikzpicture}
	\path[every node/.append style={circle, fill=black, minimum size=6pt, label distance=0pt, inner sep=0pt}]
    (0,0) node[label={-90:\(a_3\)}] (a3) {}
    (0,\radius) node[label={-90:\(a_2\)}] (a2) {}
    (0,2*\radius) node[label={90:\(a_1\)}] (a1) {}
    (2*\radius,\radius) node[label={-90:\(a_6\)}] (a6) {}
    (2*\radius,2*\radius) node[label={90:\(a_5\)}] (a5) {}
    (2*\radius,0) node[label={-90:\(a_4\)}] (a4) {}
    (-\radius,-\radius) node[red,label={-90:\(b_1\)}] (b1) {}
    (3*\radius,-\radius) node[blue,label={-90:\(c_1\)}] (c1) {};
    \draw[line width=1.2] (a1) edge (a5) edge (a6) edge[red] (b1)
    (a2) edge (a4) edge (a5) edge (a6)
    (a3) edge[red] (b1) edge (a4) edge (a6)
    (c1) edge[blue] (a4) edge[blue] (a5) edge[green!50!black] (b1);
\end{tikzpicture}
\caption{A graph describing the $[[7,1,3]]_2$ Steane code.}
\label{graphsteane}
\end{center}
\end{figure}
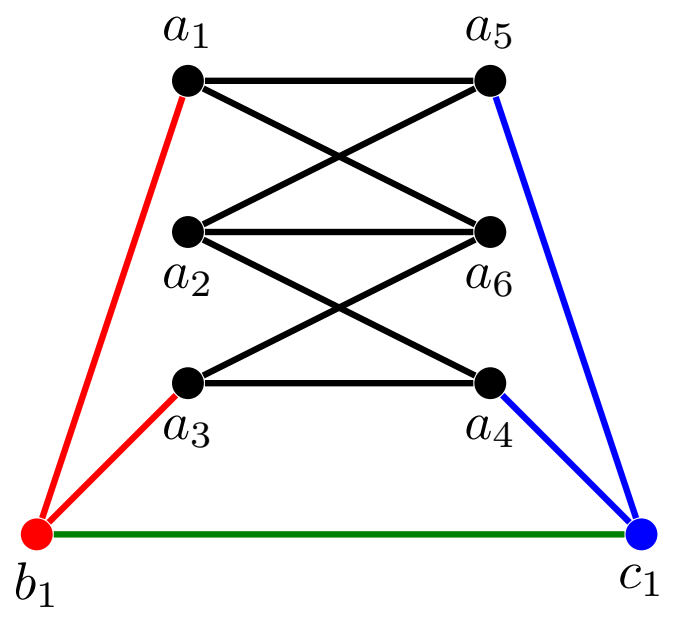
\end{example}

\subsection{Local complementation}

We define two graphs to be \emph{stabiliser-equivalent} if they give rise to equivalent stabiliser codes. Recall that two stabiliser codes are equivalent if they can be obtained from one another by
permutations of the \(n\) pairs of positions, and, for each pair, a symplectic transformation on that pair.

In the case \(k=0\) and \(p=2\), equivalent graph states have attracted attention in \cite{CP2025,localcomp,VJHBGP2025} for example.
For \(k=0\) and general primes \(p\), Bahramgiri and Beigi \cite[Theorem~5]{BB2006} describe this stabiliser-equivalence as a sequence of graph operations.
Let \(A=\left(a_{ij}\right)_{i,j\in\{1,\dots,n\}}\) denote the adjacency matrix of a graph \(G\) whose edges have weights over \(\mathbb{F}_p=\mathbb{Z}/p\mathbb{Z}\). For every vertex \(m\in\{1,\dots,n\}\) and constant \(\lambda\in\mathbb{F}_p\setminus\{0\}\), consider the following \emph{local complementation operations} on the stabiliser matrix:
    \begin{enumerate}
        \item Weighted local complementation: \[a_{ij}\mapsto a_{ij}'=a_{ij}+a_{mi}a_{mj}\quad\text{for all \(i,j\in\{1,\dots,n\}\), \(i\neq j\).}\]
        
        \item Rescaling the neighbourhood: \[a_{ij}\mapsto a_{ij}'=\begin{cases}\lambda a_{ij}&\text{if }m\in\{i,j\}\text{,}\\a_{ij}&\text{else.}\end{cases}\quad\text{for all \(i,j\in\{1,\dots,n\}\).}\]
    \end{enumerate}
They correspond to the following operations on the stabiliser code \(S\):
\begin{enumerate}
    \item 
    replacing \((v_m,v_{n+m})\) by \((v_m-v_{n+m},v_{n+m})\) within every vector \(\v\in S\),
    \item 
    replacing \((v_m,v_{n+m})\) by \((\lambda^{-1}v_m,\lambda v_{n+m})\) within every vector \(\v\in S\),
\end{enumerate}
and after each operation, putting the stabiliser matrix back in the standard form
\[G=\left(\begin{array}{c||c}I_n&A\end{array} \right)\]
through row operations, and, in the case of weighted local complementation, making the diagonal of \(A\) zero again by replacing \((v_i,v_{n+i})\) with \((v_i, v_{n+i}-a_{mi}^2v_i)\) within every vector \(\v\in S\), for all \(i\neq m\).
These two types of operations are enough to generate the entire local Clifford group:

\begin{theorem}[{\cite[Theorem~5]{BB2006}}]\label{thm:localcomplementation}
    Two (one-dimensional) graph states are equivalent if and only if there exists a sequence of local complementation operations to obtain one from the other.
\end{theorem}

\subsection{Stabiliser-equivalence for multidimensional graph states}

In this section, we work out what equivalence means for the multidimensional case. Recall from Theorem~\ref{thm:ABC} that an \([\![n,k]\!]\) stabiliser code corresponds to a graph \(\Gamma\) on \(n+k\) vertices, containing two independent sets of size \(k\) with a matching between them.
More precisely, \(\Gamma=(V,E)\) with \(V=V_A\sqcup V_B\sqcup V_C\) such that \(|V_A|=n-k\) and \(|V_B|=|V_C|=k\). The induced graphs on \(V_B\) and \(V_C\) are empty and there is a non-weighted matching between \(V_B\) and \(V_C\). This graph should be regarded as unlabelled on the vertices within each of the parts \(V_A\), \(V_B\) and \(V_C\), but this partition (colouring) is fixed. In particular, we cannot simply interchange the \(B\) and \(C\) parts.

Recall that we can always put the stabiliser matrix in the form
\[\left(\begin{array}{@{}ccc|ccc||ccc|ccc@{}}
1&\cdots&0&b_{11}&\cdots&b_{1k}&a_{11}+\sum_{i=1}^kb_{1i}c_{1i}&\cdots&a_{1r}+\sum_{i=1}^kb_{1i}c_{ri}&c_{11}&\cdots&c_{1k}\\
\vdots&\ddots&\vdots&\vdots&&\vdots&\vdots&\ddots&\vdots&\vdots&&\vdots\\
0&\cdots&1&b_{r1}&\cdots&b_{rk}&a_{r1}+\sum_{i=1}^kb_{ri}c_{1i}&\cdots&a_{rr}+\sum_{i=1}^kb_{ri}c_{ri}&c_{r1}&\cdots&c_{rk}
\end{array}\right)\]
and label the vertices of the graph as follows:
\begin{align*}
    V_A&=\{a_1,\dots,a_r\},\\
    V_B&=\{b_1,\dots,b_k\},\\
    V_C&=\{c_1,\dots,c_k\}.
\end{align*}

We consider three operations, and prove that a successive application of these three is enough to obtain every equivalent multidimensional graph state, similar to \cite[Theorem~5]{BB2006}.

Note that isomorphic graphs (but with the same \((V_A,V_B,V_C)\) partition) are equivalent. This corresponds to permutations of the pairs of positions. In what follows, isomorphic graphs are considered as the same.

\begin{theorem}\label{thm:equivoperations}
    Two multidimensional graph states are equivalent if and only if they can be obtained from one another through a sequence of the following operations.
    \begin{enumerate}
        \item (Local complementation on a vertex of \(V_B\) or \(V_C\)). Given \(b_s\in V_B\), do \[\begin{cases}
            a_{ij}\mapsto a_{ij}'=a_{ij}+b_{is}b_{js}&\text{for all \(i,j\in\{1,\dots,r\}\), \(i\neq j\),}\\
            c_{is}\mapsto c_{is}'=c_{is}-b_{is}&\text{for all \(i\in\{1,\dots,r\}\).}
        \end{cases}\]
        
        Similarly for a vertex \(c_s\in V_C\).
        
        
        In particular, we can always make a given vertex of \(V_A\) non-adjacent to a given vertex of \(V_C\) (or \(V_B\)).


        \item (Local complementation on a vertex of \(V_A\)). Given \(a_m\in V_A\), if \(c_{mj}=0\) for all \(j\in\{1,\dots,k\}\) (which can always be accomplished by operation 1), do one of the following operations:
        \begin{enumerate}
            \item Weighted complementation: \[\begin{cases}a_{ij}\mapsto a'_{ij}=a_{ij}+a_{mi}a_{mj}&\text{for all \(i,j\in\{1,\dots,r\}\), \(i\neq j\),}\\
            b_{ij}\mapsto b'_{ij}=b_{ij}+a_{mi}b_{mj}&\text{for all \(i\in\{1,\dots,r\}\) and \(j\in\{1,\dots,k\}\).}\end{cases}\]
            
            \item Rescaling the neighbourhood with \(\lambda\in\mathbb{F}_p\setminus\{0\}\): \[\begin{cases}a_{ij}\mapsto\begin{cases}\lambda a_{ij}&\text{if }m\in\{i,j\},\\a_{ij}&\text{else}\end{cases}&\text{for all \(i,j\in\{1,\dots,r\}\), \(i\neq j\)}\\ b_{ij}\mapsto\begin{cases}\lambda b_{ij}&\text{if }m=i,\\b_{ij}&\text{else}\end{cases}&\text{for all \(i\in\{1,\dots,r\}\) and \(j\in\{1,\dots,k\}\).}\end{cases}\]
        \end{enumerate}
    
        \item (Swapping a vertex of \(V_A\) with a vertex of \(V_C\)). Given \(a_m\in V_A\) and \(c_s\in V_C\), if \(c_{mj}=0\) for all \(j\in\{1,\dots,k\}\) and \(b_{ms}=1\) (which can always be accomplished by operation 1), swap \(a_m\) and \(c_s\) (only edges between \(V_A\) and \(V_C\) are affected by this swap), and change the weights of the edges between \(V_A\) and \(V_B\) as follows:
        \[b_{ij}\mapsto b_{ij}'=\begin{cases}b_{ij}&\text{if }m=i,\\-b_{ij}&\text{if }m\neq i\text{ and }s=j,\\ b_{ij}-b_{is}b_{mj}&\text{else}\end{cases}\]
        for all \(i\in\{1,\dots,r\}\) and \(j\in\{1,\dots,k\}\).
    \end{enumerate}
\end{theorem}
\begin{proof}
    We show that the first operation corresponds to adding the first column of \(B\) to the first column of \(C\), that the second operation achieves local equivalence on \(\Gamma_A\), and that the third one arises from a transposition of a pair on the first \(n-k\) positions and a pair on the last \(k\) positions.
    \begin{enumerate}
        \item Without loss of generality, we may assume that \(s=1\). Subtracting the first column of \(B\) from the first column of \(C\) gives
        \[\resizebox{\linewidth}{!}{\(\left(\begin{array}{@{}ccc|ccc||ccc|cccc@{}}
        1&\cdots&0&b_{11}&\cdots&b_{1k}&a_{11}+\sum_{i=1}^kb_{1i}c_{1i}&\cdots&a_{1r}+\sum_{i=1}^kb_{1i}c_{ri}&c_{11}-b_{11}&c_{12}&\cdots&c_{1k}\\
        \vdots&\ddots&\vdots&\vdots&&\vdots&\vdots&\ddots&\vdots&\vdots&\vdots&&\vdots\\
        0&\cdots&1&b_{r1}&\cdots&b_{rk}&a_{r1}+\sum_{i=1}^kb_{ri}c_{1i}&\cdots&a_{rr}+\sum_{i=1}^kb_{ri}c_{ri}&c_{r1}-b_{r1}&c_{r2}&\cdots&c_{rk}
        \end{array}\right).\)}\]
        Writing \(a'_{ij}=a_{ij}+b_{i1}b_{j1}\) and \(c'_{i1}=c_{i1}-b_{i1}\), this is equal to
        \[\left(\begin{array}{@{}ccc|ccc||ccc|ccc@{}}
        1&\cdots&0&b_{11}&\cdots&b_{1k}&a'_{11}+\sum_{i=1}^kb_{1i}c'_{1i}&\cdots&a_{1r}'+\sum_{i=1}^kb_{1i}c'_{ri}&c'_{11}&\cdots&c'_{1k}\\
        \vdots&\ddots&\vdots&\vdots&&\vdots&\vdots&\ddots&\vdots&\vdots&&\vdots\\
        0&\cdots&1&b_{r1}&\cdots&b_{rk}&a'_{r1}+\sum_{i=1}^kb_{ri}c'_{1i}&\cdots&a'_{rr}+\sum_{i=1}^kb_{ri}c'_{ri}&c'_{r1}&\cdots&c'_{rk}
        \end{array}\right).\]

        \item Without loss of generality, we may assume that \(m=1\). Since \(c_{1j}=0\) for all \(j\in\{1,\dots,k\}\) and \(a_{11}=0\), the stabiliser matrix is given by
        \[\resizebox{\linewidth}{!}{\(\left(\begin{array}{@{}cccc|ccc||cccc|ccc@{}}
        1&0&\cdots&0&b_{11}&\cdots&b_{1k}&0&a_{12}+\sum_{i=1}^kb_{1i}c_{2i}&\cdots&a_{1r}+\sum_{i=1}^kb_{1i}c_{ri}&0&\cdots&0\\
        0&1&\cdots&0&b_{21}&\cdots&b_{2k}&a_{21}&a_{22}+\sum_{i=1}^kb_{2i}c_{2i}&\cdots&a_{2r}+\sum_{i=1}^kb_{2i}c_{ri}&c_{21}&\cdots&c_{2r}\\
        \vdots&\vdots&\ddots&\vdots&\vdots&&\vdots&\vdots&\vdots&\ddots&\vdots&\vdots&&\vdots\\
        0&0&\cdots&1&b_{r1}&\cdots&b_{rk}&a_{r1}&a_{r2}+\sum_{i=1}^kb_{ri}c_{2i}&\cdots&a_{rr}+\sum_{i=1}^kb_{ri}c_{ri}&c_{r1}&\cdots&c_{rk}
        \end{array}\right).\)}\]

        \begin{enumerate}[(a)]
            \item 
        Subtracting column \(n+1\) from column \(1\) gives the matrix
        \[\resizebox{\linewidth}{!}{\(\left(\begin{array}{@{}cccc|ccc||cccc|ccc@{}}
        1&0&\cdots&0&b_{11}&\cdots&b_{1k}&0&a_{12}+\sum_{i=1}^kb_{1i}c_{2i}&\cdots&a_{1r}+\sum_{i=1}^kb_{1i}c_{ri}&0&\cdots&0\\
        -a_{12}&1&\cdots&0&b_{21}&\cdots&b_{2k}&a_{21}&a_{22}+\sum_{i=1}^kb_{2i}c_{2i}&\cdots&a_{2r}+\sum_{i=1}^kb_{2i}c_{ri}&c_{21}&\cdots&c_{2r}\\
        \vdots&\vdots&\ddots&\vdots&\vdots&&\vdots&\vdots&\vdots&\ddots&\vdots&\vdots&&\vdots\\
        -a_{1r}&0&\cdots&1&b_{r1}&\cdots&b_{rk}&a_{r1}&a_{r2}+\sum_{i=1}^kb_{ri}c_{2i}&\cdots&a_{rr}+\sum_{i=1}^kb_{ri}c_{ri}&c_{r1}&\cdots&c_{rk}
        \end{array}\right),\)}\]
        which, after row operations, becomes
        \[\resizebox{\linewidth}{!}{\(\left(\begin{array}{@{}cccc|ccc||cccc|ccc@{}}
        1&0&\cdots&0&b'_{11}&\cdots&b'_{1k}&a'_{11}&a'_{12}+\sum_{i=1}^kb'_{1i}c_{2i}&\cdots&a'_{1r}+\sum_{i=1}^kb'_{1i}c_{ri}&0&\cdots&0\\
        0&1&\cdots&0&b'_{21}&\cdots&b'_{2k}&a'_{21}&a'_{22}+\sum_{i=1}^kb'_{2i}c_{2i}&\cdots&a'_{2r}+\sum_{i=1}^kb'_{2i}c_{ri}&c_{21}&\cdots&c_{2r}\\
        \vdots&\vdots&\ddots&\vdots&\vdots&&\vdots&\vdots&\vdots&\ddots&\vdots&\vdots&&\vdots\\
        0&0&\cdots&1&b'_{r1}&\cdots&b'_{rk}&a'_{r1}&a'_{r2}+\sum_{i=1}^kb'_{ri}c_{2i}&\cdots&a'_{rr}+\sum_{i=1}^kb'_{ri}c_{ri}&c_{r1}&\cdots&c_{rk}
        \end{array}\right),\)}\]
        where \(a'_{ij}=a_{ij}+a_{1i}a_{1j}\) and \(b'_{ij}=b_{ij}+a_{1i}b_{1j}\). After making the diagonal entries of \(A'\) equal to zero, this is a graph state.

        \item 
        Multiplying column 1 by \(\lambda^{-1}\) and multiplying column \(n+1\) by \(\lambda\), yields the matrix
        \[\resizebox{\linewidth}{!}{\(\left(\begin{array}{@{}cccc|ccc||cccc|ccc@{}}
        \lambda^{-1}&0&\cdots&0&b_{11}&\cdots&b_{1k}&0&a_{12}+\sum_{i=1}^kb_{1i}c_{2i}&\cdots&a_{1r}+\sum_{i=1}^kb_{1i}c_{ri}&0&\cdots&0\\
        0&1&\cdots&0&b_{21}&\cdots&b_{2k}&\lambda a_{21}&a_{22}+\sum_{i=1}^kb_{2i}c_{2i}&\cdots&a_{2r}+\sum_{i=1}^kb_{2i}c_{ri}&c_{21}&\cdots&c_{2r}\\
        \vdots&\vdots&\ddots&\vdots&\vdots&&\vdots&\vdots&\vdots&\ddots&\vdots&\vdots&&\vdots\\
        0&0&\cdots&1&b_{r1}&\cdots&b_{rk}&\lambda a_{r1}&a_{r2}+\sum_{i=1}^kb_{ri}c_{2i}&\cdots&a_{rr}+\sum_{i=1}^kb_{ri}c_{ri}&c_{r1}&\cdots&c_{rk}
        \end{array}\right).\)}\]
        Multiplying the first row by \(\lambda\), we get
        \[\resizebox{\linewidth}{!}{\(\left(\begin{array}{@{}cccc|ccc||cccc|ccc@{}}
        1&0&\cdots&0&\lambda b_{11}&\cdots&\lambda b_{1k}&0&\lambda a_{12}+\sum_{i=1}^k\lambda b_{1i}c_{2i}&\cdots&\lambda a_{1r}+\sum_{i=1}^k\lambda b_{1i}c_{ri}&0&\cdots&0\\
        0&1&\cdots&0&b_{21}&\cdots&b_{2k}&\lambda a_{21}&a_{22}+\sum_{i=1}^kb_{2i}c_{2i}&\cdots&a_{2r}+\sum_{i=1}^kb_{2i}c_{ri}&c_{21}&\cdots&c_{2r}\\
        \vdots&\vdots&\ddots&\vdots&\vdots&&\vdots&\vdots&\vdots&\ddots&\vdots&\vdots&&\vdots\\
        0&0&\cdots&1&b_{r1}&\cdots&b_{rk}&\lambda a_{r1}&a_{r2}+\sum_{i=1}^kb_{ri}c_{2i}&\cdots&a_{rr}+\sum_{i=1}^kb_{ri}c_{ri}&c_{r1}&\cdots&c_{rk}
        \end{array}\right).\)}\]
        This is a graph state when substituting every \(a_{1i}=a_{i1}\) for \(\lambda a_{1i}\) and every \(b_{1j}\) for \(\lambda b_{1j}\).
        \end{enumerate}

        \item Without loss of generality, we may assume that \(m=s=1\). Since \(c_{1j}=0\) for all \(j\in\{1,\dots,k\}\) and \(b_{11}=1\) and \(a_{11}=0\), the stabiliser matrix is given by
        \[\resizebox{\linewidth}{!}{\(\left(\begin{array}{@{}cccc|cccc||cccc|ccc@{}}
        1&0&\cdots&0&1&b_{12}&\cdots&b_{1k}&0&a_{12}+\sum_{i=1}^kb_{1i}c_{2i}&\cdots&a_{1r}+\sum_{i=1}^kb_{1i}c_{ri}&0&\cdots&0\\
        0&1&\cdots&0&b_{21}&b_{22}&\cdots&b_{2k}&a_{21}&a_{22}+\sum_{i=1}^kb_{2i}c_{2i}&\cdots&a_{2r}+\sum_{i=1}^kb_{2i}c_{ri}&c_{21}&\cdots&c_{2r}\\
        \vdots&\vdots&\ddots&\vdots&\vdots&\vdots&&\vdots&\vdots&\vdots&\ddots&\vdots&\vdots&&\vdots\\
        0&0&\cdots&1&b_{r1}&b_{r2}&\cdots&b_{rk}&a_{r1}&a_{r2}+\sum_{i=1}^kb_{ri}c_{2i}&\cdots&a_{rr}+\sum_{i=1}^kb_{ri}c_{ri}&c_{r1}&\cdots&c_{rk}
        \end{array}\right).\)}\]
        Swapping column \(1\) and the first column of \(B\) and swapping column \(n+1\) and the first column of \(C\), gives
        \[\resizebox{\linewidth}{!}{\(\left(\begin{array}{@{}cccc|cccc||cccc|cccc@{}}
        1&0&\cdots&0&1&b_{12}&\cdots&b_{1k}&0&a_{12}+\sum_{i=1}^kb_{1i}c_{2i}&\cdots&a_{1r}+\sum_{i=1}^kb_{1i}c_{ri}&0&0&\cdots&0\\
        b_{21}&1&\cdots&0&0&b_{22}&\cdots&b_{2k}&c_{21}&a_{22}+\sum_{i=1}^kb_{2i}c_{2i}&\cdots&a_{2r}+\sum_{i=1}^kb_{2i}c_{ri}&a_{21}&c_{22}&\cdots&c_{2r}\\
        \vdots&\vdots&\ddots&\vdots&\vdots&\vdots&&\vdots&\vdots&\vdots&\ddots&\vdots&\vdots&\vdots&&\vdots\\
        b_{r1}&0&\cdots&1&0&b_{r2}&\cdots&b_{rk}&c_{r1}&a_{r2}+\sum_{i=1}^kb_{ri}c_{2i}&\cdots&a_{rr}+\sum_{i=1}^kb_{ri}c_{ri}&a_{r1}&c_{r2}&\cdots&c_{rk}
        \end{array}\right).\)}\]
        Bringing it in row echelon form, we get
        \[\resizebox{\linewidth}{!}{\(\left(\begin{array}{@{}cccc|cccc||cccc|cccc@{}}
        1&0&\cdots&0&1&b_{12}'&\cdots&b_{1k}'&0&c_{12}+a_{21}+\sum_{i=2}^kb'_{1i}c_{2i}&\cdots&c_{1r}+a_{r1}+\sum_{i=2}^kb'_{1i}c_{ri}&0&0&\cdots&0\\
        0&1&\cdots&0&b_{21}'&b'_{22}&\cdots&b'_{2k}&c_{21}&a_{22}+b'_{21}a_{21}+\sum_{i=2}^kb'_{2i}c_{2i}&\cdots&a_{2r}+b'_{21}a_{r1}+\sum_{i=2}^kb'_{2i}c_{ri}&a_{21}&c_{22}&\cdots&c_{2r}\\
        \vdots&\vdots&\ddots&\vdots&\vdots&\vdots&&\vdots&\vdots&\vdots&\ddots&\vdots&\vdots&\vdots&&\vdots\\
        0&0&\cdots&1&b'_{r1}&b'_{r2}&\cdots&b'_{rk}&c_{r1}&a_{r2}+b'_{r1}a_{21}+\sum_{i=2}^kb_{ri}c_{2i}&\cdots&a_{rr}+b'_{r1}a_{r1}+\sum_{i=2}^kb'_{ri}c_{ri}&a_{r1}&c_{r2}&\cdots&c_{rk}
        \end{array}\right).\)}\]
        where \(b'_{ij}\) is defined as in the statement. This is equivalent to a graph state where \(a_1\) and \(c_1\) are swapped.
\end{enumerate}

    Now, consider two equivalent stabiliser codes. By definition, they can be obtained from one another by a permutation of the coordinate positions and a symplectic transformation on each of the coordinate pairs.
    
    \begin{itemize}
        \item Permutations on the positions correspond to graph isomorphism and the second operation.
        \item Symplectic transformations on the first \(n-k\) positions correspond to the third operation, by Theorem~\ref{thm:localcomplementation}.
        \item Symplectic transformations on the last \(k\) positions can be realized by the first operation since \(\mathrm{Sp}(2,p)=\mathrm{SL}(2,p)\) is generated multiplicatively by \[\begin{pmatrix}1&1\\0&1\end{pmatrix}\text{ and }\begin{pmatrix}1&0\\1&1\end{pmatrix}.\]
    \end{itemize}
    Note that the order of the operations is important, since the first one is sometimes necessary before applying the other two.
\end{proof}

\begin{corollary}
    Two binary multidimensional graph states are equivalent if and only if they can be obtained from one another through a sequence of the following operations.
    \begin{enumerate}
        \item (Local complementation on a vertex of \(V_B\) or \(V_C\)). Given \(b_s\in V_B\), do \[\begin{cases}
            a_{ij}\mapsto a_{ij}'=a_{ij}+b_{is}b_{js}&\text{for all \(i,j\in\{1,\dots,r\}\), \(i\neq j\),}\\
            c_{is}\mapsto c_{is}'=c_{is}+b_{is}&\text{for all \(i\in\{1,\dots,r\}\).}
        \end{cases}\]
        
        Similarly for a vertex \(c_s\in V_C\).

        \item (Local complementation on a vertex of \(V_A\)). Given \(a_m\in V_A\), if \(c_{mj}=0\) for all \(j\in\{1,\dots,k\}\) (which can always be accomplished by operation 1), do a local complementation (but keeping the induced subgraph on \(V_B\) empty): \[\begin{cases}a_{ij}\mapsto a'_{ij}=a_{ij}+a_{mi}a_{mj}&\text{for all \(i,j\in\{1,\dots,r\}\), \(i\neq j\),}\\
            b_{ij}\mapsto b'_{ij}=b_{ij}+a_{mi}b_{mj}&\text{for all \(i\in\{1,\dots,r\}\) and \(j\in\{1,\dots,k\}\).}\end{cases}\]
    
        \item (Swapping a vertex of \(V_A\) with a vertex of \(V_C\)). Given \(a_m\in V_A\) and \(c_s\in V_C\), if \(c_{mj}=0\) for all \(j\in\{1,\dots,k\}\) and \(b_{ms}=1\) (which can always be accomplished by operation 1), swap \(a_m\) and \(c_s\) (only edges between \(V_A\) and \(V_C\) are affected by this swap), and change the weights of the edges between \(V_A\) and \(V_B\) as follows:
        \[b_{ij}\mapsto b_{ij}'=\begin{cases}b_{ij}&\text{if }m=i\text{ or }s=j,\\ b_{ij}+b_{is}b_{mj}&\text{else}\end{cases}\]
        for all \(i\in\{1,\dots,r\}\) and \(j\in\{1,\dots,k\}\).
    \end{enumerate}
\end{corollary}


\begin{remark}
    We cannot simply swap \(V_B\) and \(V_C\) to obtain an equivalent state. For example, if we change the roles of \(B\) and \(C\) in the \([\![4,1,2]\!]_2\) stabiliser code with stabiliser matrix
    \[\left(\begin{array}{@{}ccc|c||ccc|c@{}}
    1&0&0&0&0&0&0&1\\
    0&1&0&1&1&1&1&1\\
    0&0&1&1&1&1&1&1
    \end{array}\right),\]
    then we get the \([\![4,1,1]\!]_2\) stabiliser code
    with stabiliser matrix
    \[\left(\begin{array}{@{}ccc|c||ccc|c@{}}
    1&0&0&1&0&1&1&0\\
    0&1&0&1&0&1&1&1\\
    0&0&1&1&0&1&1&1
    \end{array}\right).\]
    These codes have a different minimum distance, so they cannot be equivalent.
    Note that \(A=O\), so \(\Gamma_A\) is edgeless.
\end{remark}

\begin{remark}\label{rem:lc}
    It is not necessarily so that \(\Gamma_A\) is determined up to local complementation (independent of \(\Gamma_B\) and \(\Gamma_C\)), as it is for the case \(k=0\). 
    For example, the \([\![5,1,2]\!]_2\) stabiliser code with stabiliser matrix
    \[\left(\begin{array}{@{}cccc|c||cccc|c@{}}
    1&0&0&0&0&0&1&1&1&1\\
    0&1&0&0&0&1&0&1&1&1\\
    0&0&1&0&1&0&0&1&0&1\\
    0&0&0&1&1&0&0&0&1&1
    \end{array}\right)\]
    corresponds to a graph where \(\Gamma_A\) is a triangle plus one pendant vertex (leaf). However, one can check (by computer; there are 24 graph states in its equivalence orbit) that it is not equivalent to a multidimensional graph state where the base graph is a path on four vertices.
\end{remark}

\section{Stabiliser codes as quantum sets of lines}

Let \(\PG(r-1,p)\) denote the \((r-1)\)-dimensional projective space over the finite field \(\FF_p\), obtained from an \(r\)-dimensional vector space over \(\FF_p\). Given a stabiliser matrix \(G\) of an \([\![n,k]\!]_p\) code, we can look at its \(2n\) columns as points in PG\((n-k-1,p)\).
These columns appear in pairs, so we can consider the (multi)set of lines\footnote{The span of two columns is not necessarily a line if they are scalar multiples. However, this means that the stabiliser code is \emph{degenerate}: deleting that coordinate gives a smaller code with an equal or larger minimum distance.} \(\left\{\ell_i\colon\,i\in\{1,\dots,n\}\right\}\), where $\ell_i$ is spanned by the $i$-th and $(i+n)$-th column of \(G\). Such a set of lines is called \emph{quantum}:

\begin{definition}\label{def:qsol}
    A (multi)set of lines \(\ell_1,\dots,\ell_n\) in \(\PG(r-1,p)\) is \emph{quantum} if there exists a stabiliser matrix \(G\) with elements in \(\FF_p\) such that its $i$th and $(i+n)$th column span \(\ell_i\) for all \(i\in\{1,\dots,n\}\).
\end{definition}

Glynn, Gulliver, Maks and Gupta \cite{GGMG} gave a description of these quantum sets of lines for $p=2$.
\begin{theorem}[{\cite[Theorem~1.13]{GGMG}, see also \cite[Theorem~3.7]{BCH2023}}]\label{thm:QSOL}
    A set of lines in PG$(r-1,2)$ is quantum if and only if every codimension two subspace is skew to an even number of its lines (see Figure~\ref{fig:QSOL}).
    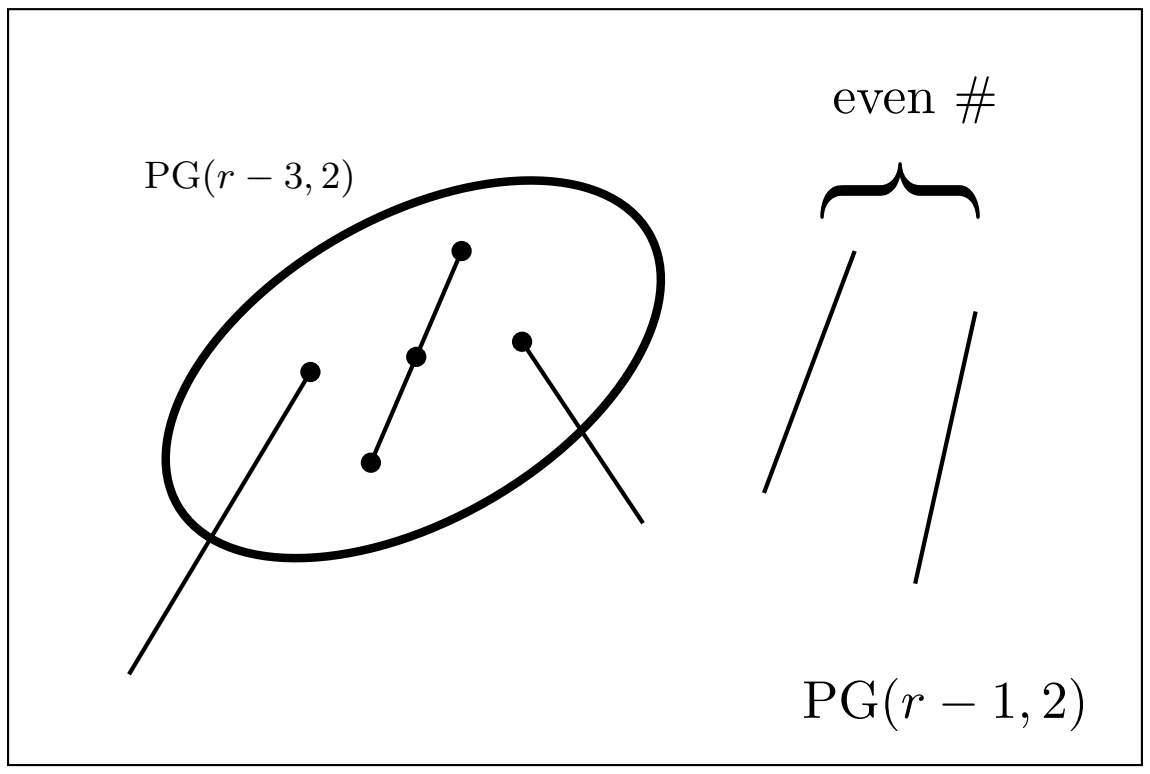
\begin{figure}[H]
    \centering
    \begin{tikzpicture}[scale=2]
\draw[rotate=120,ultra thick] (1.8,.3) ellipse (.5 and .9) {};
\draw[thick] (-2.1,.4) -- (-1.5,1.4);
\draw[thick] (-1.3,1.1) -- (-1,1.8);
\draw[thick] (-.8,1.5) -- (-.4,.9);
\draw[thick] (0,1) -- (.3,1.8);
\draw[thick] (.5,.7) -- (.7,1.6);
\draw[fill] (-1.5,1.4) {} circle (.03)
(-1.3,1.1) {} circle (.03)
(-1.15,1.45) {} circle (.03)
(-1,1.8) {} circle (.03)
(-.8,1.5) {} circle (.03);
\node[scale=.8] at (-1.7,2.04) {\small\(\PG(r-3,2)\)};
\node at (.6,.3) {\(\PG(r-1,2)\)};
\node[scale=3,rotate=-90] at (.45,2) {\(\{\)};
\node at (.5,2.3) {even \#};
\draw (-2.5,.1) rectangle (1.25,2.6);
\end{tikzpicture}
    \caption{A quantum set of lines in PG$(r-1,2)$ has the property that every codimension two subspace is skew to an even number of its lines.}
    \label{fig:QSOL}
    \end{figure}
\end{theorem}

Binary stabiliser codes are equivalent if and only if their quantum sets of lines are projectively equivalent. However, if \(p\geq3\), then projective equivalence is no longer the same as local Clifford equivalence. Indeed, multiplying a column with a nonzero scalar different from one, changes the symplectic inner product.
Also Theorem~\ref{thm:QSOL} fails when \(p\geq3\).
Still, we can use the geometry to determine the minimum distance of a code, even if \(p\geq3\). This is explained in the next section.

\subsection{The minimum distance of a quantum set of lines}

Recall that $r$ points are {\em independent} if they span an $(r-1)$-dimensional subspace; they are {\em dependent} otherwise. The following theorem is a generalisation of \cite[Theorem 3.8]{BCH2023}, which itself was a refinement of \cite[Proposition 1.14]{GGMG}.

\begin{theorem}\label{thm:QSOLd}
    Let \(X\) be a quantum set of lines coming from an \([\![n,k,d]\!]_p\) stabiliser code. If \(k=0\), then \(d\) is the minimum size of a set of dependent points on distinct lines of \(X\).
    If \(k>0\), then \(d\) is the minimum size of a set of dependent points on distinct lines of \(X\) such that there is no hyperplane containing these points and all the lines of \(X\) that do not contain one of these points.
\end{theorem}
\begin{proof}
    This is a straight-forward generalisation of the proof of \cite[Theorem 3.8]{BCH2023}, replacing \(\FF_2\) by \(\FF_p\). The main idea is that
    \[S^\perp=\left\{\v\in\FF_p^{2n}\colon\,\sum_{i=1}^n\u_iv_{n+i}-v_i\u_{n+i}=0\right\}\] where \(\u_i\) is the \(i\)th column of \(G\), together with the observation that \(\u_iv_{n+i}-v_i\u_{n+i}\in\ell_i\).
\end{proof}


\subsection{Projecting quantum sets of lines}\label{sec:projection}

The projection of an object in geometry is a common operation. In this section, we show how a quantum set of lines \(X\) projects onto a quantum set of lines \(X'\). In terms of the stabiliser matrix, a projection is equivalent to deleting rows. To see this, consider a basis in which the final \(s\) coordinates span the subspace \(\pi\) from which we are projecting. The projection is nothing more than deleting those coordinates. In the stabiliser matrix, this removes the final \(s\) rows and we get an \((n-k-s)\times2n\) matrix \(G'\). The following theorem demonstrates how we can bound the minimum distance of the projected code knowing only \(X\) and \(\pi\).

\begin{theorem}\label{thm:proj}
Let \(X\) be a quantum set of \(n\) lines in \(\PG(n-k-1,p)\) and let \(\pi\) be a subspace skew to the lines of \(X\). The projection of \(X\) from \(\pi\) is a quantum set of lines. Its minimum distance is equal to the minimum size of a set of points \(\p_1,\dots,\p_d\) on distinct lines such that they are either dependent or their span intersects \(\pi\), that is, \(\sum_{i=1}^d\p_i\in\pi\), without there being a hyperplane containing \(\p_1,\dots,\p_d\) and all the lines of \(X\) that do not contain one of these points.
\end{theorem}
\begin{proof}
    Projection is equivalent to deleting a row of a stabiliser matrix, so it preserves the quantum property.
    The claim about the minimum distance follows from Theorem~\ref{thm:QSOLd}, together with the fact that the pre-image of dependent projected points corresponds to points that sum to a point of \(\pi\) or sum to zero.
\end{proof}

We now demonstrate how the graphs associated to \(X\) and its projection \(X'\) are related.
Suppose that the stabiliser matrix is in the standard form
\[G=\left(\begin{array}{c|c||c|c}I_{n-k}&B&A+BC^T&C\end{array}\right)\]
from Theorem~\ref{thm:ABC}, and let \(B=\{\e_1,\dots,\e_r\}\) be the standard basis of \(\FF_p^r\). The elements of \(B\) are points on different lines. Let \(\{\f_1,\dots,\f_s\}\) be a basis of \(\pi\). We can extend it with elements of \(B\) to complete it to a basis of \(\FF_p^r\). Suppose, without loss of generality, that \(C=\{\e_1,\dots,\e_{r-s},\f_1,\dots,\f_s\}\) does the trick. The change of basis matrix from \(B\) to \(C\) is of the form
\[M=\left( \begin{array}{c|c}  
I_{r-s} & D_1 \\ \hline
O & D_2 \\
\end{array} \right).\]
Let
\[G=\left( \begin{array}{c|c|c||c|c|c}  
I_{r-s} & O & B_1 & A_{11}+B_1C_1^T & A_{12}+B_1C_2^T & C_1 \\ \hline
O & I_s & B_2 & A_{21}+B_2C_1^T & A_{22}+B_2C_2^T & C_2 \\
\end{array} \right)
\]
be the stabiliser matrix associated to \(X\). It is equivalent to the matrix
\[\resizebox{\linewidth}{!}{\(MG=\left( \begin{array}{c|c|c||c|c|c}  
I_{r-s} & D_1 & B_1+D_1B_2 & A_{11}+B_1C_1^T +D_1(A_{21}+B_2C_1^T) & A_{12}+B_1C_2^T+D_1(A_{22}+B_2C_2^T) & C_1+D_1C_2 \\ \hline
O & D_2 & D_2B_2 & D_2(A_{21}+B_2C_1^T) & D_2(A_{22}+B_2C_2^T) & D_2C_2 \\
\end{array} \right)
\)}\]
where the last \(s\) positions correspond to \(\pi\). Hence the matrix associated to \(X'\) is
\[\resizebox{\linewidth}{!}{\(
G'=\left( \begin{array}{c|c|c||c|c|c}  
I_{r-s} & D_1 & B_1+D_1B_2 & A_{11}+B_1C_1^T +D_1(A_{21}+B_2C_1^T) & A_{12}+B_1C_2^T+D_1(A_{22}+B_2C_2^T) & C_1+D_1C_2
\end{array} \right).
\)}\]
Comparing this with the standard form
\[
G'=\left(\begin{array}{c|c||c|c} I_{n-k} & B' & A'+B'C'^T & C' \end{array} \right)
\]
from Theorem~\ref{thm:ABC}, we have
\[A'=A_{11}-B_1C_2^TD_1^T-D_1C_2B_1^T-D_1(A_{22}+B_2C_2^T+C_2B_2^T)D_1^T,\]
\[B'=\left( \begin{array}{c|c}  
D_1 & B_1+D_1B_2
\end{array} \right)\text{ and }C'=\left( \begin{array}{c|c}  
A_{12}+B_1C_2^T+D_1(A_{22}+B_2C_2^T) & C_1+D_1C_2
\end{array}\right).\]
In particular, if we start from a one-dimensional graph state, we have
\[A'=A_{11}+D_1A_{22}D_1^T\text{, }\quad B'=D_1\text{ and }C'=A_{12}+D_1A_{22}.\]
Moreover, if the \(s\) vertices whose adjacency is given by \(A_{22}\), is an independent set, then \(A_{22}=O\) and the graph \(\Gamma'\) is obtained from \(\Gamma\) by making that independent set equal to \(V_C\) and adding the \(s\) vertices \(V_B\), the weight of the edges being determined by \(\pi\).

\subsection{Extending and shortening quantum sets of lines}\label{sec:extension}

The standard form
\[G=\left(\begin{array}{c|c||c|c}I_{n-k}&B&A+BC^T&C\end{array}\right)\]
from Theorem~\ref{thm:ABC} gives rise to a construction to extend or shorten quantum sets of lines by one line. We can add columns \(\mathbf{b}\) and \(\mathbf{c}\) to \(B\) and \(C\) respectively to obtain the stabiliser matrix
\[G'=\left(\begin{array}{c|c|c||c|c|c}I_{n-k}&B&\mathbf{b}&A+BC^T+\mathbf{b}\mathbf{c}^T&C&\mathbf{c}\end{array}\right)\]
of an \([\![n+1,k+1,d']\!]\) code. Vice versa, we can remove any two corresponding columns in \(B\) and \(C\) to shorten the code.

Let \(C\) and \(C'\) be the codes generated by \(G\) and \(G'\) respectively.

\begin{lemma}\label{lem:extend}
    Let \(\u,\v\in\FF_p^n\) and \(c\in\FF_p\). If \(\left(\begin{array}{c|c||c|c}\u&0&\v&c\end{array}\right)\in C'\), then \(\left(\begin{array}{c||c}\u&\v\end{array}\right)\in C\).
\end{lemma}
\begin{proof}
    If \(\left(\begin{array}{c|c||c|c}\u&0&\v&c\end{array}\right)\in C'\), then \(\left(\begin{array}{c|c||c|c}\u&0&\v&c\end{array}\right)=\x G'\) for some \(\x\in\FF_p^r\). Since the coordinate on position \(n+1\) is zero, we have that \(\x\mathbf{b}=0\) and hence \(\x\mathbf{b}\mathbf{c}^T=0\). Thus, \(\left(\begin{array}{c||c}\u&\v\end{array}\right)=\x G\).
\end{proof}

\begin{lemma}\label{lem:mindist}
    If \(k\geq1\), then the above extension operation can increase the minimum distance by at most one, that is, \(d(C')\leq d(C)+1\).
\end{lemma}
\begin{proof}
    First observe that the operation changes the parity check matrix of Lemma~\ref{lem:pcm1} into
    \[H'=\left(\begin{array}{c|c|c||c|c|c}
    I_{n-k}&&&A+CB^T+\mathbf{c}\b^T&&\\
    \hline
    &I_k&&C^T&&\\
    \hline
    &&1&\c^T&&\\
    \hline
    &&&-B^T&I_k&\\
    \hline
    &&&-\b^T&&1
    \end{array}\right).\]
    (Empty blocks correspond to zeroes.)
    Choose a minimum weight vector \(\left(\begin{array}{c||c}\u&\v\end{array}\right)\) that is contained in the centraliser of the original code, but not contained in the code itself (for this we need \(k\geq1\)). It is a linear combination of rows of the parity check matrix. The corresponding vector in the new parity check matrix is either \(\left(\begin{array}{c|c||c|c}\u&0&\v&0\end{array}\right)\) or \(\left(\begin{array}{c|c||c|c}\u&0&\v&0\end{array}\right)+\left(\begin{array}{c|c|c||c|c|c}\0&\0&0&\b^T&\0&0\end{array}\right)\).
    
    In the first case, the (symplectic) weight does not change. By Lemma~\ref{lem:extend}, the vector in the new parity check matrix is not contained in the code.
    
    In the second case, the vector \[\left(\begin{array}{c|c||c|c}\u&0&\v&0\end{array}\right)+\left(\begin{array}{c|c|c||c|c|c}\0&\0&0&\b^T&\0&0\end{array}\right)+\left(\begin{array}{c|c|c||c|c|c}\0&\0&0&-\b^T&\0&1\end{array}\right)=\left(\begin{array}{c|c||c|c}\u&0&\v&1\end{array}\right)\] belongs to the centraliser and has weight one more than \(\left(\begin{array}{c||c}\u&\v\end{array}\right)\). By Lemma~\ref{lem:extend}, it is not contained in the code.
\end{proof}

The extension operation can be described geometrically if \(p=2\) (this condition is needed for step 3 below). Given a quantum set of \(n\) lines in \(\PG(r-1,2)\), where \(r=n-k\), we can add a line and adjust the other lines such that they form a quantum set of lines again, as follows:

\begin{enumerate}
    \item Choose a basis \(\{\mathbf{e}_1,\dots,\mathbf{e}_r\}\) of \(r\) points on \(r\) different lines \(\ell_1,\dots,\ell_r\).
    \item Choose two different points \(\b\) and \(\c\) such that \(\forall i\in\{1,\dots,r\}\colon\,\b\in\ell_i\setminus\{\e_i\}\implies c_i=0\). 
    Let \(\ell_{n+1}=\langle\b,\c\rangle\) be the line spanned by them.
    \item For each \(i\in\{1,\dots,r\}\) replace the line \(\ell_i=\langle\e_i,\a_i\rangle\) by the line \(\ell_i'=\langle\e_i,\a_i+c_i\b\rangle\).
    \item Now \(\{\ell_1',\dots,\ell_r'\}\cup\{\ell_{r+1},\dots,\ell_n\}\cup\{\ell_{n+1}\}\) is a quantum set of \(n+1\) lines in \(\PG(r-1,2)\).
\end{enumerate}

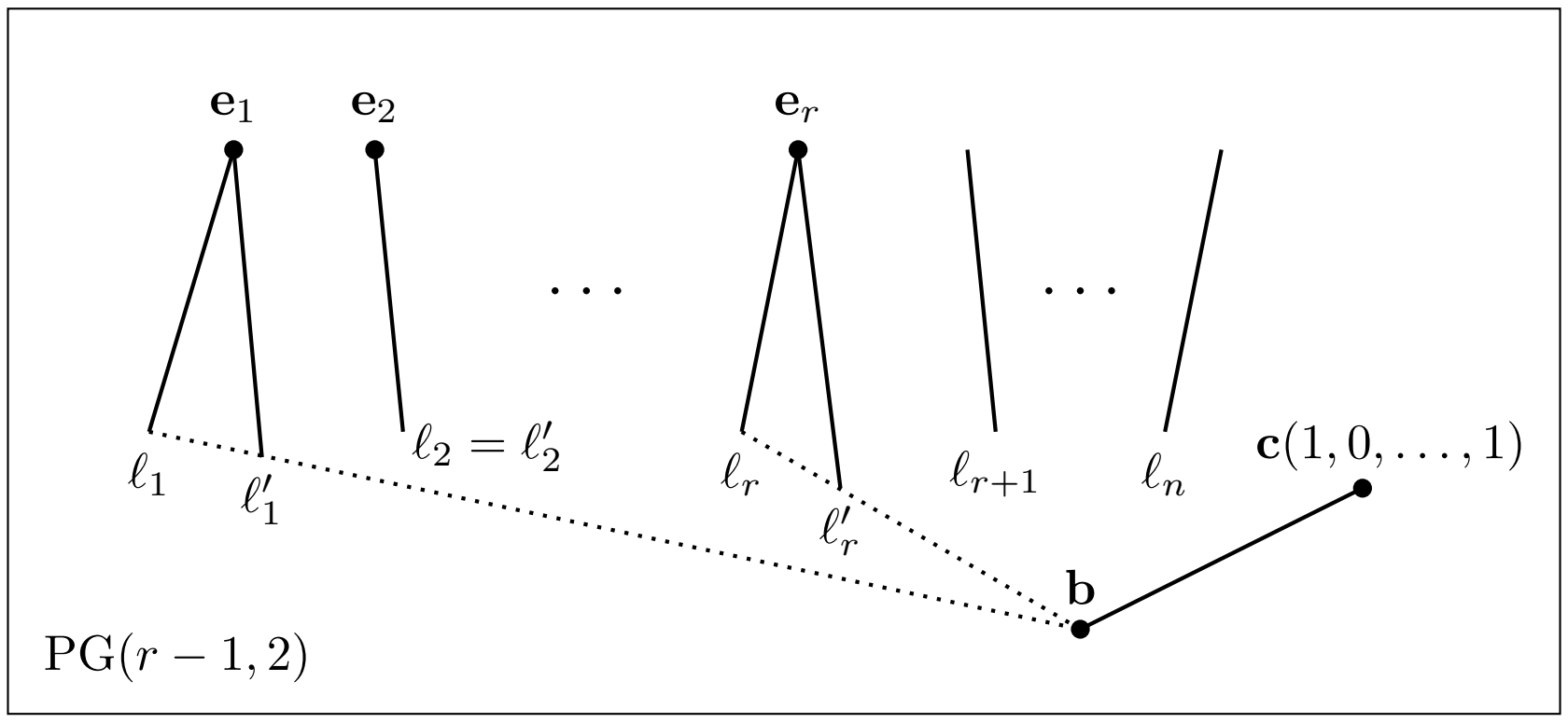
\begin{figure}[H]
    \centering
    \begin{tikzpicture}[scale=2]
\draw[thick] (0,1) -- (-.3,0);
\draw[thick] (.5,1) -- (.6,0);
\draw[thick] (2,1) -- (1.8,0);
\draw[thick] (2.6,1) -- (2.7,0);
\draw[thick] (3.5,1) -- (3.3,0);
\draw[thick] (4,-.2) -- (3,-.7);
\draw[thick,dotted] (-.3,0) -- (3,-.7);
\draw[thick,dotted] (1.8,0) -- (3,-.7);
\draw[thick] (0,1) -- (0.1,-.09);
\draw[thick] (2,1) -- (2.15,-.2);
\draw[fill] (0,1) {} circle (.03)
(.5,1) {} circle (.03)
(2,1) {} circle (.03)
(4,-.2) {} circle (.03)
(3,-.7) {} circle (.03);
\node at (0,1.15) {\(\mathbf{e}_1\)};
\node at (-.3,-.15) {\(\ell_1\)};
\node at (0.1,-.24) {\(\ell_1'\)};
\node at (0.5,1.15) {\(\mathbf{e}_2\)};
\node at (.9,-.05) {\(\ell_2=\ell_2'\)};
\node at (1.25,.5) {\Large\(\dots\)};
\node at (2,1.15) {\(\mathbf{e}_r\)};
\node at (1.8,-.15) {\(\ell_r\)};
\node at (2.15,-.35) {\(\ell_r'\)};
\node at (2.7,-.15) {\(\ell_{r+1}\)};
\node at (3.3,-.15) {\(\ell_n\)};
\node at (3,.5) {\Large\(\dots\)};
\node at (4.1,-.05) {\(\mathbf{c}(1,0,\dots,1)\)};
\node at (3,-.55) {\(\mathbf{b}\)};
\node at (-.2,-.8) {\(\PG(r-1,2)\)};
\draw (-.8,-1) rectangle (4.7,1.5);
\end{tikzpicture}
\end{figure}

Conversely, we can also remove a line from a quantum set of \(n\) lines in \(\PG(r-1,2)\), where \(r=n-k\) and \(k\geq1\), such that it is still a quantum set of lines, as follows:
\begin{enumerate}
    \item Choose a basis \(\{\mathbf{e}_1,\dots,\mathbf{e}_r\}\) of \(r\) points on \(r\) different lines \(\ell_1,\dots,\ell_r\).
    \item Pick another line \(\ell\notin\{\ell_1,\dots,\ell_r\}\) from the quantum set of lines and two different points \(\b\) and \(\c\) on that line \(\ell\) such that \(\forall i\in\{1,\dots,r\}\colon\,\b\in\ell_i\setminus\{\e_i\}\implies c_i=0\).
    \item For each \(i\in\{1,\dots,r\}\) replace the line \(\ell_i=\langle\e_i,\a_i\rangle\) by the line \(\ell_i'=\langle\e_i,\a_i+c_i\b\rangle\).
    \item Now \(\{\ell_1',\dots,\ell_r'\}\cup\{\ell_{r+1},\dots,\ell_n\}\setminus\{\ell\}\) is a quantum set of \(n-1\) lines in \(\PG(r-1,2)\).
\end{enumerate}

Note that, in the above operation, \(\mathbf{b}\) and \(\mathbf{c}\) not interchangeable. The following example illustrates this.

\begin{example}
    Consider the star graph \(K_{1,4}\). It corresponds to the graph state with stabiliser matrix
    \[\left(\begin{array}{@{}ccccc||ccccc@{}}
    1&0&0&0&0&0&1&1&1&1\\
    0&1&0&0&0&1&0&0&0&0\\
    0&0&1&0&0&1&0&0&0&0\\
    0&0&0&1&0&1&0&0&0&0\\
    0&0&0&0&1&1&0&0&0&0
    \end{array}\right)\]
    and the corresponding set of lines is a point pencil through \(\e_1\), whose projection from \(\e_1\) is a frame in the hyperplane \(X_1=0\). Suppose that we want to extend the set with another line through \(\e_1\), say the line \(\ell=\{\e_1,\b,\c\}\) with \(\b=(0,1,1,0,0)\) and \(\c=(1,1,1,0,0)\). If we follow the above operation, we get the stabiliser matrix
    \[\left(\begin{array}{@{}ccccc|c||ccccc|c@{}}
    1&0&0&0&0&0&0&1&1&1&1&1\\
    0&1&0&0&0&1&0&1&1&0&0&1\\
    0&0&1&0&0&1&0&1&1&0&0&1\\
    0&0&0&1&0&0&1&0&0&0&0&0\\
    0&0&0&0&1&0&1&0&0&0&0&0
    \end{array}\right)\]
    but if we would change the roles of \(\b\) and \(\c\) in the construction, we would get the matrix
    \[\left(\begin{array}{@{}ccccc|c||ccccc|c@{}}
    1&0&0&0&0&1&0&0&0&1&1&0\\
    0&1&0&0&0&1&1&1&1&0&0&1\\
    0&0&1&0&0&1&1&1&1&0&0&1\\
    0&0&0&1&0&0&1&0&0&0&0&0\\
    0&0&0&0&1&0&1&0&0&0&0&0
    \end{array}\right).\]
    These quantum sets of lines are not equivalent because the second one contains the line \(\{\e_2,\e_3\}\) with multiplicity \(2\), while the first one contains every line only once.
\end{example}

Note that the choice of the basis \(\{\e_1,\dots,\e_r\}\) matters. See also Remark~\ref{rem:lc}.

\section{Application: classifying stabiliser codes up to equivalence}\label{sec:application}

Using Theorem~\ref{thm:ABC} and Theorem~\ref{thm:equivoperations}, we can classify all multidimensional graph states of a certain small size up to equivalence by computer.

We provide tables containing the number of graph states, up to equivalence, for \(k\in\{0,1\}\) and \(p\in\{2,3\}\), generalising computational results for \(k=0\) and \(p=2\) by Danielsen and Parker \cite{classification} and for \(k=0\) and \(p=3\) by Danielsen \cite{classification1}. These tables were obtained with SageMath \cite{sagemath} by a straight-forward implementation of the local complementation strategy from Theorem~\ref{thm:equivoperations}. Our implementation can most probably be optimised with better code. But we do not focus on performance. Our goal is mostly to illustrate that, with the newly established graph formalism, it becomes easier to classify stabiliser codes with a given set of parameters for larger dimensions (\(k\geq1\)) as well.

\begin{table}[H]
    \centering
    \begin{tabular}{ccccccccccccc}
        \hline
        \(d\backslash n\) &1&2&3&4&5&6&7&8&9&10&11&12\\
        \hline
        1 &1&1&2&3&6&11&26&59&182&675&3990&45,144\\
        2 &&1&1&3&4&13&29&107&416&2618&27,445&615,180\\
        3 &&&&&1&1&4&11&69&577&11,202&467,519\\
        4 &&&&&&1&&5&8&120&2506&195,456\\
        5 &&&&&&&&&&&1&63\\
        6 &&&&&&&&&&&&1\\
        \hline
        All & 1&2&3&6&11&26&59&182&675&3990&45,144&1,323,363\\
        \hline
    \end{tabular}
    \caption{The number of \([\![n,0,d]\!]_2\) stabiliser codes up to equivalence. Table taken from \cite[Table~4]{classification}.}
    \label{tab:countk=0p=2}
\end{table}

\begin{table}[H]
    \centering
    \begin{tabular}{ccccccccc}
        \hline
        \(d\backslash n\) &1&2&3&4&5&6&7&8\\
        \hline
        1 & 1&2&5&11&29&78&260&1023\\
        2 & &&&2&6&35&169&1170\\
        3 & &&&&1&2&19&178\\
        \hline
        All &1&2&5&13&36&115&448&2371\\
        \hline
    \end{tabular}
    \caption{The number of \([\![n,1,d]\!]_2\) stabiliser codes up to equivalence.}
    \label{tab:countk=1p=2}
\end{table}

\begin{table}[H]
    \centering
    \begin{tabular}{ccccccc}
        \hline
        \(d\backslash n\) &1&2&3&4&5&6\\
        \hline
        1 &1&1&2&3&7&13\\
        2 &&1&1&3&5&20\\
        3 &&&&1&1&5\\
        4 &&&&&&1\\
        \hline
        All & 1&2&3&7&13&39\\
        \hline
    \end{tabular}
    \caption{The number of \([\![n,0,d]\!]_3\) stabiliser codes up to equivalence. The last row agrees with the values in \cite[Table~2]{classification1}.}
    \label{tab:countk=0p=3}
\end{table}

\begin{table}[H]
    \centering
    \begin{tabular}{cccccc}
        \hline
        \(d\backslash n\) &1&2&3&4&5\\
        \hline
        1 & 1&2&5&12&37\\
        2 & &&1&4&20\\
        3 & &&&&1\\
        \hline
        All &1&2&6&16&58\\
        \hline
    \end{tabular}
    \caption{The number of \([\![n,1,d]\!]_3\) stabiliser codes up to equivalence.}
    \label{tab:countk=1p=3}
\end{table}

\section{Final comments and further work}

We proved that an $[\![n,k,d]\!]_p$ stabiliser code is equivalent to a graph on $n+k$ vertices, which has two independent sets on $k$ vertices with a matching in between, and whose edges are weighted with elements of ${\mathbb F}_p$ (Theorem~\ref{thm:ABC}). We translated local Clifford equivalence to a set of operations on such graphs (Theorem~\ref{thm:equivoperations}).

It does not appear to be easy to deduce the minimum distance of a stabiliser code directly from the associated graph. However, it can be calculated using the geometry of quantum sets of lines (Theorem~\ref{thm:QSOLd}).

The graph formalism allowed us to define two operations on quantum sets of lines to obtain other sets of quantum sets of lines. As an application, we classified stabiliser codes for some given small parameters.

We propose the following directions for future work:

\begin{enumerate}
    \item Can one extend the ideas of this article to prime powers that are not primes? For geometrical objects equivalent to stabiliser codes of fields of non-prime order, see \cite{BMS2025,Heinrich2021}.


    \item Extend the tables in Section~\ref{sec:application}. Something that might prove interesting, is an algorithm to calculate a canonical representative of a stabiliser code, similarly to the work of Feulner \cite{Feulner} for linear codes.

    \item Can we lower bound the minimum distance of a stabiliser code using graph parameters of the corresponding graph state? What are families of graphs for which the associated stabiliser codes have a large minimum distance?

    \item The projection and extension operations from Section~\ref{sec:projection} and Section~\ref{sec:extension} can be used to find codes with a large minimum distance in a recursive way. For example, we can follow the path
    \[[\![1,0,1]\!]_2\,\longrightarrow\,[\![2,1,1]\!]_2\,\longrightarrow\,[\![2,0,2]\!]_2\,\longrightarrow\dots\longrightarrow\,[\![n,1,d]\!]_2\,\longrightarrow\,[\![n,0,d']\!]_2\]
    while optimizing for the largest possible minimum distance.
    
    To our knowledge, the existence of a \([\![14,3,5]\!]_2\) stabiliser code is still open. We applied the extension operation on every \([\![11,0,5]\!]_2\) stabiliser code to see if a \([\![14,3,5]\!]_2\) could be obtained by walking the path
    \[[\![11,0,5]\!]_2\quad\longrightarrow\quad[\![12,1,5]\!]_2\quad\longrightarrow\quad[\![13,2,4]\!]_2\quad\longrightarrow\quad[\![14,3,5]\!]_2.\]
    After going through all possibilities, we conclude that it is impossible. This gives some evidence (but not a proof) that a \([\![14,3,5]\!]_2\) stabiliser code does not exist.
    
    The drawback with this strategy is that there are many ways to choose a line that extends a given quantum set of lines. For the given values, this was possible, but in general it seems like walking every possible path towards a \([\![14,3,5]\!]_2\) stabiliser code is infeasible on a normal computer.
\end{enumerate}

\paragraph{Acknowledgements}
We would like to thank Oriol Farr\`as and Cristina Fern\'andez for their careful reading of the doctoral thesis of Ricard Vilar, on which the results in this paper are partially based.

\paragraph{AI declaration}
We did not use any AI-assisted technologies in the research, writing, or editing of the manuscript.


\vfill
\noindent
\begin{minipage}[t]{0.5\textwidth}
\textsc{Simeon Ball}\\
\textsc{\small Department of Mathematics}\\[-1mm]
\textsc{\small Universitat Politècnica de Catalunya}\\[-1mm]
\textsc{\small C. Pau Gargallo 14, 08028 Barcelona, Spain}\\[-1mm]
{\it\small E-mail:} {\tt\small simeon.michael.ball@upc.edu}\\[2mm]

\textsc{Robin Simoens}\\
\textsc{\small Department of Mathematics: Analysis, Logic}\\[-1mm]
\textsc{\small and Discrete Mathematics}\\[-1mm]
\textsc{\small Ghent University}\\[-1mm]
\textsc{\small Krijgslaan 281, 9000 Gent, Belgium}\\[-1mm]
\textsc{\small Department of Mathematics}\\[-1mm]
\textsc{\small Universitat Politècnica de Catalunya}\\[-1mm]
\textsc{\small C. Pau Gargallo 14, 08028 Barcelona, Spain}\\[-1mm]
{\it\small E-mail:} {\tt\small Robin.Simoens@UGent.be}\\[2mm]
\end{minipage}
\hfill
\begin{minipage}[t]{0.47\textwidth}
\textsc{G\"ulizar G\"unay}\\
\textsc{\small Institute of Analysis and Number Theory}\\[-1mm]
\textsc{\small Graz University of Technology}\\[-1mm]
\textsc{\small M\"unzgrabenstra{\ss}e 36, 8010 Graz, Austria}\\[-1mm]
{\it\small E-mail:} {\tt\small gunaymert@math.tugraz.at}\\[2mm]

\textsc{Ricard Vilar}\\
\textsc{\small Department of Mathematics}\\[-1mm]
\textsc{\small Universitat Politècnica de Catalunya}\\[-1mm]
\textsc{\small C. Pau Gargallo 14, 08028 Barcelona, Spain}\\[-1mm]
{\it\small E-mail:} {\tt\small ricard.vilar@upc.edu}
\end{minipage}

\end{document}